\documentclass[10pt]{article}
 
\usepackage{latexsym}
\usepackage{amsmath}
\usepackage{amsfonts}
\usepackage{eufrak}
\usepackage{overcite}
\usepackage{theorem}
\theorembodyfont{\upshape}

\newtheorem{theorem}{Theorem}

\newtheorem{lemma}{Lemma}

\def\Bbb#1{{\mathbb#1}}

\def\RRd{\Bbb{R}^d}

\def\ZZd{\Bbb{Z}^d}

\def\G{{\cal G}}
\def\H{{\cal H}}

\def\half{\frac{1}{2}}

\def\Ltwo{L^2(\mu_{\rho})}
\def\xr{\chi(\rho)}

\def\murho{{\mu}_{\rho}}
\def\gradk{\nabla_{e_k}\xi(0)}

\newcommand{\rint}[2]{\langle #1 , #2 \rangle_{\rho}}

\newcommand{\dpair}[2]{{\langle\langle #1 , #2 \rangle\rangle}_\rho}
\newcommand{\dnorm}[1]{{\langle\langle #1 \rangle\rangle}_{\rho}}

\begin{document}


\begin{center}
\vspace*{5mm}

{\bf \Large On the Symmetry of the Diffusion Coefficient in Asymmetric Simple Exclusion. }\\
\vspace{6mm}
Michail Loulakis\footnote[1]{Statistical Laboratory, Centre for Mathematical Sciences, Wilberforce Road, Cambridge CB3 0WB, United Kingdom. E-mail: {\tt michail@statslab.cam.ac.uk}}

\end{center}
\vspace{6mm}
\small
\hspace{1cm} \begin{minipage}[t]{4.9in}
ABSTRACT: We prove the symmetry of the diffusion coefficient that appears in the fluctuation-dissipation theorem for asymmetric simple exclusion processes.
\end{minipage}
\\
\normalsize
\vspace{8mm}

\noindent
\section{Introduction}

The fluctuation-dissipation theorem has been an instrumental tool in the analysis of non-gradient interacting particle systems. In the context of asymmetric simple exclusion processes (ASEP) it has been used to establish the hydrodynamic limit for the mean-zero ASEP {\cite{X}}, as well as the diffusive incompressible limit \cite{EMY}, the first order corrections to the hydrodynamic limit \cite{LOY1}, the equilibrium fluctuations of the density field \cite{CLO}, and the diffusive hydrodynamic limit when the initial density profile is constant along the direction of the drift \cite{LSV}, for the general ASEP in $d\ge 3$.
\\
\\
The fluctuation-dissipation theorem consists of decomposing the (normalized) particle currents $(w_i)_{1\le i\le d}$ (which are not of gradient form) into gradients of the occupation variable and a rapidly fluctuating term. With a suitable interpretation it can be formulated in the following equation:
\begin{equation}
w_i=\sum_{j=1}^d D_{ij}(\eta(0)-\eta(e_j))+Lu_i.
\label{flux}
\end{equation}

\noindent
The matrix $D=(D_{ij})_{1\le i,j\le d}$ is called the  {\em diffusion coefficient} and it naturally appears in the PDEs that arise in the hydrodynamic limit. Explicit and variational formulae for $D$ are available \cite{LOY2} and it is known to be a smooth function of the particle density \cite{LOV}. In this article we prove that $D$ is symmetric, thus answering the question raised by Landim, Olla and Yau in ref. \citen{LOY1} and \citen{LOY2}.


\section{Notation and Results}

Let us fix a finite range probability measure $p(\cdot)$ on $\mathbb{Z}^d$, with $p(0)=0$. We denote by $L$ the generator of the simple exclusion process associated to $p(\cdot)$.  $L$  acts on local functions on the state space ${\Bbb X}=\{0,1\}^{\ZZd}$ according to: 

\begin{equation}
Lf(\xi)=\sum_{x,y}p(y-x) \,\xi(x) \,(1-\xi(y))\,(f(\xi^{x,y})-f(\xi)) 
\label{generator}
\end{equation}
\noindent
where: \hspace{20mm}
 $\xi^{x,y}(z) = \begin{cases}
                 \xi(z)& \text{if $z\neq x,y$}\/,\\
                 \xi(x)& \text{if $z=y$},\\
                 \xi(y)& \text{if $z=x$}.
                \end{cases} $
\\
\\
The symmetric and the anti-symmetric part of $p(\cdot)$ will be denoted by $a(\cdot)$ and $b(\cdot)$ respectively:
\[
a(x)=\frac{p(x)+p(-x)}{2} , \qquad b(x)=\frac{p(x)-p(-x)}{2}.
\] 
In order to avoid degeneracies we will assume that the random walk in
$\ZZd$ with one step transition probabilities $a(y-x)$ is irreducible, i.e. $\{x: a(x)>0\}$ generates the group $\ZZd$.
 An equivalent formulation of this assumption is that the matrix $S=(S_{ij})_{1\le i,j\le d}$ defined by $S_{ij}=\half\sum p(z)z_iz_j$ is invertible.
\\
\\
The symmetric part of the generator (denoted by $L^s$) is given by (\ref{generator}) with $p(\cdot)$ replaced by $a(\cdot)$. The measures $\murho \, (0\! \leq \rho \!\leq 1),$ defined as Bernoulli
products of parameter $\rho$ over the sites of $\ZZd$ are invariant under the dynamics. We will denote expectations under $\murho$ by $\langle\cdot\rangle_\rho$ and inner products in $\Ltwo$ by $\rint{\cdot}{\cdot}$ .
\\
\\
The adjoint of $L$ in $\Ltwo$ is the generator $L^{*}$ of the simple exclusion process associated to the law $p^{*}(x)=p(-x)$. Local functions form a core of both $L$ and $L^*$, and thus $L^s$ extends to a self-adjoint operator in $\Ltwo$.
\\
\\
The particle current along the direction $e_i$ is given by: 
\begin{equation}
W_i=\half\sum_z p(z)z_i\xi(0)(1-\xi(z))-p(-z)z_i\xi(z)(1-\xi(0)).
\label{current}
\end{equation}

\noindent
Equation (\ref{flux}) is to be understood in the Hilbert space of fluctuations, which we define next. Let $\G_\rho$ be the space of local functions $g$ such that:
\[
\langle g \rangle_\rho=0 \qquad {\text{, and}} \qquad \left.\frac{d}{d\theta}\langle g  \rangle_\theta\right|_{\theta=\rho}=0.
\]
For a $g\in\G_{\rho}$ we define $\tau_x g=g(\tau_x\xi)$, where $\tau_x \xi(z)=\xi(x+z)$. For any $f\in\G_\rho$ and $i\in\{1,\ldots,d\}$ we define:
\[
 \langle g, f\rangle_{\rho,0}:=\sum_x\rint{g}{\tau_xf}, \qquad t_i(g)=\rint{g}{\sum_x x_i\xi(x)}.
\]
 Set $\xr=\rho(1-\rho)$ and define  
\begin{equation}
\dnorm{g}=\sup_{\alpha\in\RRd}\left(2\sum_{i=1}^d\alpha_it_i(g)-\xr\alpha\cdot S\alpha\right)+\sup_{f\in\G_\rho}\left(2\langle g,f\rangle_{\rho,0}-\langle f, (-L^s)f\rangle_{\rho,0}\right),
\label{var}
\end{equation}
\noindent

The Hilbert space of fluctuations $\H(\rho)$ is defined as the closure of $\G_\rho$ under $\dnorm{\cdot}^{1/2}$. If we denote by $\H_0$ the space generated by gradients of the occupation variable: $\H_0=\{\sum\alpha_i(\xi(e_i)-\xi(0));\ \alpha\in\RRd\}$, then by Theorem 5.9 in ref. \citen{EMY} we have:
\begin{equation}
\H(\rho)=\overline{\H_0+L\G_\rho}.
\label{decomposition}
\end{equation}
\noindent
Notice that unless $\sum zp(z)=0$ the currents $W_i$ do not belong to the space $\G_\rho$. Therefore we define the {\em normalized} currents $w_i\in\G_\rho$ by:
\begin{equation}
w_i=W_i-\langle W_i \rangle_{\rho}- (\xi(0)-\rho)\left.\frac{d}{d\theta}\langle W_i\rangle_{\theta}\right|_{\theta=\rho}.
\label{normalised}
\end{equation}
\noindent
According to (\ref{decomposition}) there exist coefficients $(D_{ij})_{1\le i,j\le d}$ (which depend on $\rho$) such that:
\begin{equation}
w_i-\sum_{j=1}^d D_{ij}\times(\xi(0)-\xi(e_j)) \in \overline{L\G_\rho}.
\label{difcoef}
\end{equation} 
\noindent
The matrix $D=(D_{ij})_{1\le i,j\le d}$ is called the diffusion coefficient of the simple exclusion process. Landim, Olla and Yau proved explicit and variational formulae for $D$ in ref. \citen{LOY2}. In the same paper, as well as in ref. \citen{LOY1}, the authors question whether there exists a choice of $p(\cdot)$ such that $D$ is asymmetric. The result of this article is the following theorem:
\begin{theorem}
The diffusion coefficient $D$ defined in (\ref{difcoef}) is always a symmetric matrix.
\label{the}
\end{theorem} 

\section{Some Properties of $\H(\rho)$}
In this section we review some properties of the Hilbert space of fluctuations that will be useful in the proof of Theorem \ref{the}. We begin with the following lemma.
\begin{lemma}
If $g\in\G_\rho$ and $h\in\ZZd$, then $\tau_hg=g$ in $\H(\rho)$.
\label{translation}
\end{lemma}  
{\bf Proof:} In view of (\ref{var}) it suffices to show that:
\[
(i)\ \langle \tau_h g-g, f\rangle_{\rho,0}=0, \forall f\in\G_\rho, \qquad
(ii)\ t_i(\tau_h g-g)=0,  i=1,\ldots,d.
\]
\noindent
Using the translation invariance of $\murho$ property $(i)$ follows immediately, while 
\[
\sum_x \rint{\tau_h g-g}{x_i\xi(x)}=h_i\sum_x\rint{g}{\xi(x)}.
\]
\noindent
The last expression is trivially zero if $\rho\in\{0,1\}$, while otherwise by differentiating with respect to $\theta$ both sides of the following identity 
\[
\langle g\rangle_\theta=\int g(\xi)\prod_{x\in {\text{supp}}(g)} \left(\frac{\theta}{\rho}\right)^{\xi(x)}\left(\frac{1-\theta}{1-\rho}\right)^{1-\xi(x)}\, d\murho(\xi),
\]
\noindent
we get 
\begin{equation}
\sum_x\rint{g}{\xi(x)}=\rho\langle g\rangle_{\rho}+\xr \left. \frac{d}{d\theta}\langle g\rangle_{\theta}\right|_{\theta=\rho}=0, 
\label{degree1}
\end{equation}
thus establishing $(ii).\qquad \Box$
\\
\\
The following lemma is a generalisation of (5.1) in ref. \citen{LOY2} to the general ASEP, and can be proved by polarization of (\ref{var}). The details are left to the reader.
\begin{lemma}
Let $g,f\in\G_\rho$ and set $\gradk=\xi(0)-\xi(e_k)$ for $k=1,\ldots, d$. Then:
\begin{eqnarray*}
(i)& & \dpair{\gradk}{\nabla_{e_\ell}\xi(0)}=\xr(S^{-1})_{k\ell},\\
(ii)& &\dpair{\gradk}{Lg}=-\dpair{\gradk}{L^{*}g}=\sum_{\ell=1}^d(S^{-1})_{k\ell}\langle w_\ell,g\rangle_{\rho,0},\\
(iii)& &\dpair{\gradk}{L^sg}=0,\\
(iv)& &\dpair{L^sg}{f}=-\langle g,f\rangle_{\rho,0}.
\end{eqnarray*}
\label{products}
\end{lemma}
\section{The Diffusion Matrix}
Recall the definition of the normalised currents $w_i$ given in  (\ref{current}) and (\ref{normalised}). It follows by elementary algebra and Lemma \ref{translation} that $w_i=W_i^s-h_i$,
where 
\[
W_i^s(\xi)=\half\sum z_ia(z)(\xi(0)-\xi(z))
\]
is the current of the symmetric simple exclusion with generator $L^s$, and 
\[
h_i(\xi)=\sum_z z_ib(z)(\xi(0)-\rho)(\xi(z)-\rho).
\]
Hence, the normalized currents for the reversed process are given by $w_i^{*}=W_i^s+h_i$. Let now $C(\rho)$ (resp. $C_{*}(\rho)$) be the real vector space generated by the currents $\{w_i; i=1,\ldots,d\}$ (resp. $\{w_i^{*}; i=1,\ldots,d\}$).
\\
\\
We define the linear operator $T$ (resp. $T^{*}$) on $C(\rho)+L\G_\rho$ (resp. $C_{*}(\rho)+L^{*}\G\rho$) by:
\[
T(\sum_{i=1}^{d}\alpha_iw_i+Lg)=\sum_{i,k=1}^d\alpha_iS_{ik}\gradk+L^sg,
\]
\[
T^{*}(\sum_{i=1}^{d}\alpha_iw^{*}_i+L^{*}g)=\sum_{i,k=1}^d\alpha_iS_{ik}\gradk+L^sg.
\]
By Theorem 5.9 in ref. \citen{EMY} we have: $\H(\rho)=\overline{C(\rho)+L\G_\rho}=\overline{C_{*}(\rho)+L^{*}\G\rho}$. Now, just as in Lemma 5.4 in ref. \citen{LOY2}, $T$ and $T^{*}$ are norm bounded by 1, hence they can be extended to $\H(\rho)$.  Furthermore, it follows easily by computations based on Lemma \ref{products} that $T^{*}$ is the adjoint of $T$ with respect to $\dpair{\cdot}{\cdot}$ and $T^{*}\gradk$ is orthogonal to $\overline{L\G_\rho}$. Hence by (\ref{difcoef}):
\[
\dpair{w_i}{T^{*}\gradk}=\sum_{i=1}^{d}D_{ij}\dpair{\nabla_{e_j}\xi(0)}{T^{*}\gradk},
\]
and thus by Lemma \ref{products}(i):
\[
\xr I_d= D\cdot Q
\]
where the matrix $Q=(Q_{jk})_{1\le j,k\le d}$ is given by: 
\[
Q_{jk}=\dpair{T\nabla_{e_j}\xi(0)}{\gradk}.
\]
We are now ready to proceed with the proof of Theorem \ref{the}. 
\\
\\
{\bf Proof} (of Theorem \ref{the}): Let us denote the reflection operator on ${\Bbb X}$ by 
\[
R\xi(z)=\xi(-z).
\]
The action of $R$ is naturally extended to functions as $Rf(\xi)=f(R\xi)$. Clearly, $R^2=1$. Furthermore, the following commutation relation can be readily verified:
\begin{equation}
RL=L^{*}R.
\label{commutation}
\end{equation}
In particular $R$ commutes with $L^s$ and hence, $R$ preserves inner products in $\H(\rho)$.
\\
\\
Notice that $W_i^s$ are anti-symmetric under $R$, while $h_i$ are $R$-symmetric. Thus, 
\begin{equation}
Rw_i(\xi)=-W_i^s(\xi)-h_i(\xi)=-w_i^{*}(\xi).
\label{currev}
\end{equation}
It is a direct consequence of (\ref{commutation}), (\ref{currev}), and the observation that $R\gradk=-\gradk$ in $\H(\rho)$ that
\[
RT=T^*R.
\]
Therefore,
\begin{eqnarray*}
Q_{jk}&=&\dpair{RT\nabla_{e_j}\xi(0)}{R\gradk}\\
&=&\dpair{T^*R\nabla_{e_j}\xi_(0)}{R\gradk}\\
&=&\dpair{\nabla_{e_j}\xi(0)}{T\gradk}\\
&=&Q_{kj}.
\end{eqnarray*}
So $Q$ and thus the diffusion coefficient $D$ are symmetric matrices.\qquad $\Box$
\\
\\
{\bf Acknowledgments:} This research has been supported by a Marie Curie Fellowship  of the European Community Programme ``Improving Human Potential'' under the contract number HPMF-CT-2002-01610.


\begin{thebibliography}{7}

\bibitem[1]{CLO} Chung, C.C., Landim, C., Olla, S. : {\em Equilibrium Fluctuations of Asymmetric Simple Exclusion Processes in Dimension $d\ge 3$.} Probab. Th. Rel. Fields {\bf 119} (2001), 381--409.

\bibitem[2]{EMY} Esposito, R., Marra, R., Yau H.T. : {\em Diffusive Limit of Asymmetric Simple Exclusion.} Rev. Math. Phys. {\bf 6} (1994), 1233--1267.

\bibitem[3]{LOY1} Landim, C., Olla, S., Yau, H.T. : {\em First Order Correction for the Hydrodynamic Limit of Asymmetric Simple Exclusion Processes in Dimension $d\ge 3$.} Comm.  Pure  Appl. Math. {\bf 50} (1997), 149--203.

\bibitem[4]{LOY2} Landim, C., Olla, S., Yau, H.T. : {\em Some Properties of the Diffusion Coefficient for Asymmetric Simple Exclusion Processes.} Ann. Probab. {\bf 24} (1996), 1779--1807.

\bibitem[5]{LOV} Landim, C., Olla, S., Varadhan, S.R.S. : {\em On Viscosity and Fluctuation-Dissipation in Exclusion Processes.}, J. Stat. Phys. {\bf 115} (2004), no1-2, 323--363.  

\bibitem[6]{LSV} Landim, C., Sued, M., Valle, G. : {\em Hydrodynamic Limit of Asymmetric Exclusion Process Under Diffusive Scaling in $d\ge 3$}, Commun. Math. Phys. {\bf 249} (2004), 215--247.
 
\bibitem[7]{X} Xu, L. : {\em Diffusive Limit for a Lattice Gas with Short Range Interactions.} Ph.D. Thesis, New York University.

\end{thebibliography}
\end{document}